\documentclass[12pt,reqno]{amsart}
\usepackage[dvips]{graphics}
\usepackage{amssymb}
\usepackage{dsfont}
\usepackage{hyperref}

\numberwithin{equation}{section}

\newtheorem{thm}{Theorem}[section]
\newtheorem*{thm*}{Theorem}
\newtheorem*{thmmain*}{MAIN THEOREM}
\newtheorem{lem}[thm]{Lemma}
\newtheorem{cor}[thm]{Corollary}
\newtheorem{prop}[thm]{Proposition}
\newtheorem*{prop*}{Proposition}
\theoremstyle{definition}
\newtheorem{defn}{Definition}[section]
\theoremstyle{remark}
\newtheorem{rem}{Remark}[section]
\newtheorem{ex}[rem]{Example}

\newcommand{\tref}[1]{Theorem~\ref{#1}}
\newcommand{\cref}[1]{Corollary~\ref{#1}}
\newcommand{\pref}[1]{Proposition~\ref{#1}}
\newcommand{\dref}[1]{Definition~\ref{#1}}

\newcommand{\lref}[1]{Lemma~\ref{#1}}

\newcommand{\rref}[1]{Remark~\ref{#1}}

\def\diam{\mathop{\text{diam}}}
\def\rad{\mathop{\text{rad}}}
\def\R{\mathds{R}}
\def\dim{\mathop{\text{dim}}}

\begin{document}
\title{Building-like spaces}
\author{Andreas Balser and Alexander Lytchak}


\begin{abstract}
We study convex subsets of buildings, discuss some structural features
and derive several characterizations of buildings.
\end{abstract}

\maketitle
\renewcommand{\theequation}{\arabic{section}.\arabic{equation}}
\pagenumbering{arabic}

\section{Introduction}
The main purpose of this paper is to investigate the structure 
of convex subsets of spherical buildings. Such a convex
subset $X$ inherits from the ambient spherical building $G$ the  
following fundamental property (see Subsection \ref{basicbuild}):

$(*)$  For each $x\in X$ there is some (conicality radius)
$r_x >0$, such that for all 
$y,z\in X$ with $d(x,y)< r_x$ the triangle $xyz$ is spherical.

In order to study the  geometry of $X$ it seems natural
to forget about  the ambient space $G$ and to work directly inside $X$.
Moreover it is reasonable to investigate the slightly
bigger (synthetically defined) class of all CAT(1) spaces  
that  satisfy the above property  $(*)$ and have finite geometric 
dimension. We call such spaces \emph{building-like}.

\begin{rem}
The property $(*)$ can be regarded as a variant of a constant curvature
 $1$ condition lying between sphericality (i.e. between convex subsets
of Hilbert spheres) and local conicality (i.e. 
spherical complexes, see Subsection \ref{locallyConicalSubsection}).
Observe that local conicality has almost no implications on the topology 
(compare \cite{Berest}), in contrast to the very special topology of 
building-like spaces, see \tref{introThm2}. 
\end{rem}

The first result describes the global topology of building-like spaces, characterizes buildings among them and provides a synthetic approach
to buildings.
\begin{thm}\label{introThm2}
For a  building-like space $X$ of dimension $n$,  the following are
equivalent:

\begin{enumerate}
\item $X$ is a building;

\item $X$ is geodesically complete;

\item Each point has an antipode;

\item $X$ contains an  $n$-dimensional Euclidean sphere;

\item $X$ is not contractible;

\item The radius of $X$ is equal to  $\pi$.
\end{enumerate}
\end{thm}

Moreover we show that if an $n$-dimensional building-like space 
$X$ contains an $(n-1)$ dimensional Euclidean sphere then 
$X$ is either a building or it has radius $\frac \pi 2$.

 Now we turn to the local geometry of building-like spaces 
and prove the following structural results 
(see Section \ref{difsec}  for more details).
Each building-like space $X$ has the same local dimension at all 
points, and a (topologically and
metrically) naturally defined (thick) decomposition in cells.
The most important feature of the building-like space~$X$ is
the boundary   $\partial X$, that can be described 
(in analogy to convex subsets of Riemannian manifolds) by the property  that
it is the largest subset of $X$  whose complement is convex and everywhere
 dense. The boundary is empty iff $X$ is a building, moreover the boundary
 has the    following local description that is 
a differential analog of \tref{introThm2}: 

\begin{thm} \label{boundarythm}
Let $X$ be a building-like space, $x\in X$ a point. The following
are equivalent:

\begin{enumerate}
\item $x\in \partial X$;

\item Some geodesic  terminates in $x$;

\item There are small neighborhoods $U$ of $x$, such that
 $U \setminus \{ x \} $ is contractible;

\item The link $S_x$ is not a building;

\item The link $S_x$ is contractible.

\end{enumerate}
\end{thm}

\begin{rem}
In fact if $X$ is a convex subset of a building $G$ and not an abstract
building-like space, then  all links are  building-like 
(compare Subsection \ref{boundsubsec}) and one can
deduce \tref{boundarythm}
 from \tref{introThm2}. In this case $x\in \partial X$ iff
$\partial (S_x)$ is not empty.
\end{rem}

We finish the investigation of the local structure of building-like spaces 
by showing that spherical subsets of $X$ can be extended up to the 
boundary $\partial X$.
The corresponding result for spherical buildings is shown in
\cite[Prop.\ 3.9.1]{kleinerLeeb}. 
However, our argument (unlike the proof in \cite{kleinerLeeb}) remains
valid in the Euclidean and the hyperbolic situation (see \rref{hyperbolic}). 
The corresponding  result in the case of Euclidean or hyperbolic buildings is
probably known, but we could not find a reference. 

\begin{prop} \label{extend}
Let $X$ be a building-like space of dimension $n$. Let $C\subset X$
be a convex spherical subset. Then $C$ is contained in some $n$-dimensional
spherical subset $\tilde C$, whose boundary $\partial \tilde C$  
(as a spherical set) is contained in $\partial X$.
\end{prop}

%The equivalence $(2)\to (1)$ from \tref{introThm2} now implies:
%\begin{cor}
% A building-like space has no boundary iff it is a building.
%\end{cor}

Now we discuss the local-global equivalence of our notion.  
The main result of \cite{charneyLytchak} shows that a CAT(1)
space that is locally isometric to a building is itself a building. 
This result has the following natural extension in the setting of
 building-like spaces:

\begin{thm}   \label{LocalGlobal}
Let $X$ be a connected CAT(1) space  of dimension at least 2. 
If each point has a convex
neighborhood which is (isometric to) a building-like space,
then $X$ is building-like.  
\end{thm}

The proof of \cite[Thm.\ 4.1]{charneyLytchak} provides
 the following extension 
of the well known theorem of Cartan, saying that  the universal covering
of a complete manifold with sectional curvature $1$ is a 
Euclidean sphere.
% (The corresponding result 
%is clear in the $CAT(0)$ setting by the 
%generalized theorem of Cartan-Hadamard, see \cite{bishop}):

\begin{cor} \label{Cartan}
 Let $X$ be a complete geodesic space, such that each point $x\in X$
has a building-like neighborhood of dimension $\geq 3$. Then the 
universal covering $\tilde X$ of $X$ is CAT(1), and therefore 
 building-like.
\end{cor}

In the proof of   \tref{boundarythm}  we use  a result (in which 
 we leave the universe of building-like spaces), that we consider to be
 of independent  interest.

\begin{thm} \label{append}
Let $X$ be an $n$-dimensional CAT(1) space that has at least one pair of
 antipodes.
 Assume that each pair of
points $x,y \in X$ with $d(x,y) \geq \pi$
is contained in an $n$-dimensional sphere $S^n$.
If $X$ contains an open relatively compact subset $U$ then $X$
is a building.
\end{thm}

\begin{rem} In dimension $n=1$   the above result coincides with
 Theorem 1.1.  of  
 \cite{naganodim}.
\end{rem}

From \tref{append} we derive an unpublished result of Kleiner:

\begin{cor} \label{kleiner}
Let $X$ be a locally compact CAT(0)-space  of dimension~$n$. If each pair 
of points is contained in a flat $\R^n$, then  $X$ is a building. 
\end{cor}

\begin{rem} \label{hyperbolic}
We would like  to emphasize that most of the results discussed above
 can be word by word transferred
 from   the CAT(1) to the CAT(0) and to the CAT$(-1)$ setting.
For the corresponding notion
of Euclidean (resp. hyperbolic) building-like spaces, that 
include convex subsets
of locally conical Euclidean (resp. hyperbolic) buildings 
(in contrast to spherical building-like
spaces discussed above) all the local  results  mentioned above hold true
(except \tref{introThm2}). 
The boundary of such a (Euclidean or hyperbolic)
building-like space is empty iff the space is geodesically complete. 
In the Euclidean case (but certainly not in the hyperbolic)
this is enough to deduce 
that such a space is a Euclidean building. 
 The  corresponding statements of 
\tref{LocalGlobal} resp. \cref{Cartan} are  valid in all
dimensions in the Euclidean and hyperbolic cases.
\end{rem}

%\begin{rem}
% Convex subsets of buildings defined combinatorially have been studied
% in \cite{}. However, the combinatorially defined convexity of subsets
%  is much more  restrictive. It is equivalent to the requirement 
%  that the subset is convex (in the metric sense) and that it
%  is a union of Weyl chambers of the building. Moreover we would like
%  to mention that a similiar syntethic approach to buildings
%  from the combinatorial point of view was undertaken in \cite{}.
%\end{rem}

Our investigations were
 mainly motivated  by the question if a group that operates on a building
 by isometries and fixes some non-trivial convex subset must have a fixed point.
% In the case of combinatorially convex subsets this is known as the center
% conjecture (\cite{}).
 We refer to \cite{part1} for an answer in small dimensions and to 
 \cite{Groups} for a complete answer to the closely  related question
 about groups operating by isometries on symmetric spaces or 
Euclidean buildings.

\medbreak

Now we describe the structure of the paper: In Section  \ref{sphbas}
basic properties of spherical subsets of CAT(1) spaces are discussed 
that may also be of independent use. These results are used in Section
\ref{difsec}, the heart of this paper, where the local structure of 
building-like spaces is discussed in detail and \tref{boundarythm} is shown.
In  Section \ref{globalsec} we present the proof of \tref{LocalGlobal}.
Section \ref{beginnbuild} and Section \ref{veryend}
are independent of the rest and contain the proofs of \tref{introThm2} 
and \tref{append}. The reader only interested in these results may skip the 
rest.

\smallbreak

We would like to thank Juan Souto for useful comments on a previous version.

\section{Preliminaries}
\subsection{Notations} 
By $\R ^n$ resp. $S^n$ we  denote the Euclidean
 space resp. the Euclidean sphere of dimension $n$.
By $d$ we denote distances in metric spaces, by $B_r (x)$ we
will denote the closed ball of radius $r$ around~$x$.
Geodesics are always  parametrized by  arclength. 
For a point $x$ in a metric space $X$ we set
 $\rad _x (X) := \sup _{z\in X} d(x,z)$. The radius
of $X$ is defined by  $\rad (X) := \inf _{x\in X} \rad _x (X)$.

By $X*Z$  we denote  the spherical join of $X$ and $Z$
(\cite[pp.\ 63f]{bridsonHaefliger}).

\subsection{CAT(1) spaces}
A complete metric 
space will be called CAT(1) if each pair of 
points  with distance  $<\pi$ is connected by a geodesic and 
all triangles of perimeter less than $2\pi$ are not thicker than in
$S^2$. We refer to 
\cite[ch.\ II]{bridsonHaefliger}.  A CAT(1) space is
geodesically complete 
if each geodesic can be prolonged to an infinite local geodesic. A subset
$C$ of a CAT(1) space is convex (more precisely $\pi$-convex) 
if all points  in $C$ with distance $<\pi$
are joined by a geodesic in $C$.

 In a CAT(1) space $X$ we will denote by $S_x =S_xX$ the link at the point $x$.
For each point $x\in X$ there is a natural $1$-Lipschitz (logarithmic)
map  $p_x:X\to S_x *S^0$, where $x$ is sent to a pole of $S^0$
and the distances to $x$ are preserved. We refer to \cite{lytchakSphB}.

By $\dim (X)$ we denote the geometric dimension of $X$ studied in 
\cite{kleinerDimension}.
%The following easy lemma is shown in \cite[L.\ 11.7]{lytchakSphB} by induction:
The easy proof by induction of the following %easy 
lemma can be found in \cite{part1}:
\begin{lem} \label{easyLemma}
Let $X$ be an $n$-dimensional CAT(1) space, and let $S \subset X$ be an 
embedded $S^n$. Then for each $x\in X$ there is an antipode $y\in S$, i.e.
a point satisfying $d(x,y) \geq \pi$. Therefore we have $\rad (X) \geq \pi$.
\end{lem}

%\begin{proof}
%If $\dim S=0$, the claim is clear. 
%Choose an  arbitrary $y\in S$. Let $v \in S_y$ be the starting direction
% of $yx$. By induction, we can
%choose an 
%antipode~$w$ of $v$ in $S^{n-2} 
%=S_yS\subset S_y X$. Now we obtain an antipode of~$x$ by extending
%$xy$  inside $S$ in the
%direction of $w$.
%\end{proof}

This result implies that if in an $n$-dimensional CAT(1) space $X$ 
each point is contained in some $S^n$, then $X$ is geodesically complete.

\subsection{Buildings}
We refer to \cite[sect.\ 3]{kleinerLeeb} for an account on spherical
buildings. 
%We point out   that buildings need not be
%irreducible or thick. 
In the proofs below we will use some characterizations of buildings among finite 
dimensional geodesically
complete CAT(1) spaces derived in
 \cite{lytchakSphB} and in  \cite{charneyLytchak}.

\section{Spherical parts of CAT(1) spaces} \label{sphbas}
\subsection{Spherical subsets}  We  call a subset $T$ of a 
CAT(1) space~$X$
\emph{spherical} 
if $T$ admits an isometric embedding into a Hilbertsphere. 
A convex subset $T$ of a CAT(1) space 
$X$ is spherical iff for every choice of three points
$x_1,x_2,x_3 \in T$, the triangle $x_1x_2x_3$ is spherical. 

\begin{lem}\label{sphericalParts}
Let $X$ be a CAT(1) space. Assume that $T_1,T_2$ are convex 
spherical subsets of $X$. If for all $z_1,z_2,z_3 \in T_1 \cup T_2$ 
the triangle $z_1z_2z_3$ is spherical, then the convex hull of
$T_1 \cup T_2$ is spherical too.
\end{lem}

\begin{proof}
The last observation implies
that it is enough to prove the result in the case where $T_1$ and
 $T_2$ are of dimension $\leq 2$. One sees easily that there is an isometric
 embedding $I :T_1 \cup T_2 \to S^5$. Set $Y_i = I(T_i)$. Due to
 \cite[Thm.\ A]{lang-sch} the inverse $f=I^{-1}: (Y_1 \cup Y_2 ) \to X$
 has a $1$-Lipschitz 
 extension to the convex hull $Y$ of $Y_1 \cup Y_2$.  Remark that $f$ maps
 geodesics connecting points in $Y_1 \cup Y_2$ isometrically onto their
 images and we only  have to show that 
$f:Y\to f(Y)$ is an isometry. This reduces 
 the statement to the case where $T_1$ and $T_2$ are $1$-dimensional.
 In this case the proof can be finished in the same way as in
 \cite[L.\ 4.1]{lytchakSphB}, 
 where the special case $d(x_1,x_2)=\frac \pi 2$ for
 all $x_i\in T_i$ is covered.

In fact one shows using Toponogov, that $f$ preserves the distances to $T_i$
and that its differential is an isometry at the points of $T_i$. This 
immediately implies the result.
\end{proof}

\subsection{Locally conical spaces} \label{locallyConicalSubsection}
The following definition from \cite{charneyLytchak} is a generalization
of the concept  of a simplicial complex from \cite{ballmannSimplicial}:

\begin{defn} \label{defcon}
Let $X$ be a CAT(1) space. We  call $X$ \emph{locally conical} if for each 
$x\in X$ there is an $r_x > 0$, such that for all $y,z \in B_{r_x}(x)$ 
the triangle $xyz$ is spherical. The maximal $r_x >0$ with this property
is called the \emph{conicality radius} at $x$.
\end{defn}

\begin{rem}
 A simplicial space in the sense of \cite{ballmannSimplicial} is
locally conical. It is possible to prove the converse if $X$ is 
geodesically complete.
\end{rem}

 The above condition is equivalent to 
the requirement that $B_{r_x} (x)$ is canonically (via $p_x$) isometric
to a convex subset of $\{ x\} * S_x$.
Observe that a closed convex subset of a locally conical space is 
locally conical.

\subsection{Very spherical subsets}
 The following property was studied in \cite{lytchakSphB} under 
the additional assumption of geodesic completeness:
\begin{defn}
We say that  points $x,y$ in a CAT(1) space $X$ are \emph{close} if
for each $z\in X$ the triangle $xyz$ is spherical.  We say that
a convex subset $T \subset X$ is \emph{very spherical} in $X$ if each
pair of points of $T$ 
are close.
\end{defn}

Remark that a convex subset $T \subset X$ is very spherical 
iff for all $z\in X$ it is mapped isometrically
under the logarithmic map $p_z :X \to S_z X * S^0$.
A  very spherical subset is spherical.
For  very spherical subsets $T_1, T_2$ of $X$,  the 
convex hull of $T_1$ and $T_2$ in $X$ is a spherical subset by 
\lref{sphericalParts}.
The closure of a very spherical subset is very spherical 
and the union of a chain of very spherical 
subsets  is very spherical too. Hence every very spherical subset is contained
in a maximal one. 
 
\begin{lem}\label{closeness1}
Let $T$ be a finite-dimensional very spherical subset of $X$. Let $m$ be an inner
point of the spherical convex set $T$. If $q \in X$ is close to $m$, then
the convex hull of $q$ and $T$ is very spherical. 
\end{lem}

\begin{proof}
The convex hull $C$ of $q$ and $T$ is spherical by Lemma
\ref{sphericalParts}. 
Let~$z$ be arbitrary.  Under the map $p_z$ the subset $T$ is mapped
(isometrically) 
onto a spherical subset $\bar T \subset S_z * S^0$. Set $\bar m :=p_z (m)$
and $\bar q :=p_z (q)$.  We have $d(q,m)= d(\bar q, \bar m )$ since
$q$ and $m$ are close. Let 
$x_1\in T$ be arbitrary. Choose $x_2\in T$ such that $m$ 
is an inner point of $x_1x_2$. 
Set $\bar x _i := p_z (x_i)$.
Now in the triangle $\bar x_1 \bar x_2 \bar q$ we have
$d(\bar x_i, \bar q ) \leq d(x,q)$,
$d(\bar{x}_{1},\bar{x}_{2})=d(x_{1},x_{2})$ and $d(\bar q, \bar m)
=d(q,m)$. 
Since the triangle $x_1x_2q$ is spherical, we obtain from the CAT(1)
property in $ S_z * S^0$, that $d(\bar x_i, \bar q ) = d(x_i, q)$
(and the triangle $\bar x_1 \bar x_2 \bar q$ must be spherical).
%Using Lemma \ref{sphericalParts},  we see that $C$ is mapped under
%$p_z$ isometrically  
%onto  its image. 
Since $C$ is spherical and $x_{1}$ was arbitrary, this shows that $C$
is mapped isometrically onto its image by $p_{z}$. Since $z$ was
arbitrary, this shows that $C$ is very spherical.
\end{proof}

The proof of the following lemma is provided by the fact that a triangle
$xyz$ is spherical iff for some $m$ on $xy$ the triangles $xmz$ and $ymz$ 
are spherical and $\angle_{m}( x,z) + \angle_{m} (y,z) = \pi$.

\begin{lem}\label{closeness2}
Let $xmy$ be a geodesic in $X$. Assume that
 $m$ is close to~$x$ and to $y$. Then 
$x$ is close to $y$ iff $S_m X$ splits as $S^0 * Y$, where the sphere 
$S^0$ consists of the starting directions of $mx$ and $my$.
\qed
\end{lem}

\section{Building-like spaces} \label{beginnbuild}
\subsection{Basics} \label{basicbuild}
We recall the basic definition from the introduction:

\begin{defn} \label{defmain}
Let $X$ be a CAT(1) space. We will say that $X$ is
\emph{building-like} if it has  
finite dimension and each point $x$ has a neighborhood $B_{r_x} (x)$
consisting of points close to $x$, i.e. for all $y\in B_{r_x} (x)$ and
each $z\in X$ the triangle $xyz$ is spherical.
\end{defn}

First of all we observe that spherical buildings are building-like.
Namely for each point $x$ in a building $G$ there is some $r_x >0$
such that for each $y\in B_{r_x} (x)$ the points $x$ and $y$ are contained in
some Weyl chamber of the building. Now for each other point $z\in G$
this chamber and $z$ are contained in some apartment, hence the triangle
$xyz$ is spherical.

Observe now that the class of building-like spaces is  stable with 
respect to spherical joins. Much  more important is that a closed 
convex subset of a building-like space is building-like, in particular
convex subsets of buildings are building-like.

The definition implies that a building-like space is locally conical. 
From \lref{closeness2} we see that in a building-like space $X$
the maximal~$r_x$ satisfying the condition of \dref{defmain} coincides
with the conicality radius at $x$ from \dref{defcon} (see also
\lref{closeness3} below).

\begin{ex}  \label{nullDimCat1}
Directly from the definition we see that a $0$-dimensional CAT(1)
space is always building-like. 
It is a building iff it has at least two points. 
\end{ex}

\begin{ex} \label{oneDimCat1}
Using \tref{introThm2} one easily derives the following characterization of 
$1$-dimensional building-like spaces: Let $X$ be a CAT(1) space
of dimension $1$. Then $X$ is building-like if and only if $X$ is a building
or $X$ is a locally conical metric tree of diameter $\leq \pi$.
In the latter case we have $\rad (X) = \frac {\diam (X)} 2 \leq \frac
\pi 2$.  
\end{ex}

Lemma \ref{closeness2}  directly  implies the next

\begin{lem}\label{closeness3}
Let $\gamma = xy$ be a geodesic in a building-like space $X$,
 such that for each inner point $m$ on  
$\gamma$ the 
link $S_m$ splits as $S_m =S^{0}*Z_{m}$, with 
$S^0= \{\gamma_{m}^{+},\gamma_{m}^{-} \} $.
 Then $x$ is close to $y$. \qed
\end{lem}

 We finish the basics with some remarks about antipodes. Namely from
\dref{defmain} we deduce that if $x$ and $z$ in a building-like space
$X$ are antipodes (i.e. satisfy $d(x,z)\geq \pi$), then $d(x,z)=\pi$
and for each $y\in B_{r_x} (x)$ the broken line $xyz$ is in fact a geodesic.
In particular the set $Ant(x)$ of all antipodes of $x$ is discrete
and if it is not empty, then~$X$ contains an isometric copy of $S^0 * S_x$.

\subsection{Characterization of buildings} We turn to  \tref{introThm2}:

\begin{proof}[Proof of \tref{introThm2}]
The step $(4) \to (3)$ is given by \lref{easyLemma}.

 Assume $(3)$. Let $\gamma : [-t,0] \to X$ be a geodesic with 
$\gamma (0) =x$. Set $y= \gamma (-\varepsilon )$ for some $\varepsilon <r_x$.
Choose an antipode $z$ of $y$. There must be a geodesic from $y$ to $z$
starting in the direction of $\gamma$ and since $\gamma$ cannot branch between
$y$ and $x$, we see that $yxz$ is a geodesic. Therefore $\gamma$ does not
terminate in $x$ and $X$ is geodesically complete.

The implication $(2) \to (1)$ follows from 
\cite[Prop.\ 4.5]{lytchakSphB}, which says
that a finite dimensional geodesically complete space $X$
must be a building if  the  set of antipodes of each point
$x\in X$ is discrete.

 $(1)\to (4)$ is  clear.

The implication $(1)\to (5)$ is well-known, since the homology
of an $n$-dimensional
spherical building is non-trivial  in dimension $n$. On the
other hand if a point $x \in X$ has no antipodes, then the contraction along
geodesics starting at $x$ shows that $X$ is contractible, hence $(5)$ implies
$(3)$.

Finally $(3)\to (6)$ is clear and $(6)\to (3)$ is shown 
in \cref{smallRadius}. 
\end{proof}

\subsection{Type bounds} \label{boundsubsec}
 It seems to be difficult to distinguish between properties 
of abstract building-like spaces and convex subsets of buildings.
The only advantages of the existence 
of an ambient building we found
 is stability under ultralimits and
 building-likeness of the links. We make this more precise:

\begin{defn}
Let $W$ be a finite Coxeter group. We call a space $X$ building-like of
type bounded by $W$ if $X$ admits an isometric embedding onto a convex
subset of a building $G$ of type $W$.   
\end{defn}

If the type of $X$ is bounded by $W$ then so is the type of each  convex
subset and of each link $S_x$ of $X$. Moreover if $X_i$ is a sequence of
building-like spaces of type bounded by $W$, then so is the ultralimit
$\lim _{\omega} X_i$.

\section{Local-Global Equivalence}  \label{globalsec}
Before we are going to study the local  structure of 
building-like spaces in detail we prove \tref{LocalGlobal}, showing that the 
global property of being building-like is in fact a local one:

\begin{proof}[Proof of \tref{LocalGlobal}]
By local conicality and Lemma \ref{closeness3}  it is
enough to show that 
for arbitrary points $x,z \in X$ and a geodesic $\eta$ starting at~$z$
the triangle $xyz$ is spherical if $y= \eta (\varepsilon )$ 
and $\varepsilon >0$ is small enough. 

%A connectedness argument then shows that the triangle is spherical as
%long as $\varepsilon <r_{z}$ (the conicality radius at $z$).

First, we assume that $d(x,z) <\pi$.
The case $d(x,z)\geq \pi$ can be easily deduced afterwards using
the fact that the links are connected due to the 
assumption $\dim (X)\geq 2$ and \tref{boundarythm}.

 We can cover the geodesic  $\gamma = xz$
 by finitely many convex open building-like subsets.
Choosing $\varepsilon$ small enough, we may assume that the geodesic
$\bar \gamma =xy$ is contained
in the union of these subsets. More precisely we can find points
$x=x_0,x_1,\dotsc ,x_n=z$ on $xz$ and $x=y_0,y_1,\dotsc ,y_n =y$ on $xy$
such that for all $i$ the points $x_i,x_{i+1},y_i,y_{i+1}$ are
contained in some open convex building-like subset $U_i$.

Moving $x_{i}$ along $\gamma$, we may assume that  for each $x_i$ the
directions $\gamma ^+$ and 
$\gamma ^-$ define a factor of $S_{x_i}$ for $0<i<n$ (by local
conicality, the points where this is not true are discrete in $\gamma$). 
 Making $\varepsilon $ smaller
 % (if necessary)
  we may assume  that
$y_i$ is close to $x_i$ inside of $U_i$.  Moving $y_i$ on $\bar \gamma$
 we may assume that $ \bar \gamma ^+$ and
$\bar \gamma ^-$   define a factor of $S_{y_i}$ for $0<i<n$.

%   \lref{sphericalParts} implies
% that the convex hull $x_ix_{i+1}y_iy_{i+1}$
% is spherical. 

The triangle $xx_{1}y_{1}$ is spherical since $x_{1},y_{1}$ are
close. Assume by induction that $xx_{i}y_{i}$ is spherical. 
Our assumption on the link at $y_{i}$ implies that we can glue the
spherical triangle $x_{i}y_{i}y_{i+1}$ to obtain a spherical triangle
$xx_{i}y_{i+1}$. Similarly, we can glue the spherical triangle
$x_{i}x_{i+1}y_{i+1}$, so the triangle $xx_{i+1}y_{i+1}$ is
spherical. 
% 
% Using the fact that all $x_i$ resp. all $y_i$ are
% on a geodesic $\gamma$ resp. $\bar \gamma$ and the above assumptions
% on the links, we obtain by induction 
% on $i$, that the quadrangles $x_ix_{i+1}y_iy_{i+1}$ are 
% two-dimensional and that the triangle $x_0x_iy_i$  is spherical.
 For $i+1=n$ we get the result.
%
%Consider points $x,z$ such that there is a geodesic $\gamma$ of length
%$\pi$ connecting them (with $\gamma (0)=x,\gamma (\pi )=z$). Let $v$
%be the incoming direction of $\gamma$ in $S_{z}$, and $w\in S_{z}$
%such that $\angle (v,w)<\frac{\pi}{2}$ and a geodesic $\bar{\gamma}$
%starts at $z$ in the direction of $w$. Let $y$ be a point on
%$\bar{\gamma}$ with $d(y,z)<r_{z}$. Above, we have shown that the triangle
%$yz\gamma (\varepsilon )$ is spherical for every $\varepsilon >0$.
%Letting $\varepsilon \rightarrow 0$, we find that there is a geodesic
%of length $\pi$ connecting $x,z$ with initial direction $w$ at $z$. 
%Since our assumptions imply that $S_{z}$ is connected, this shows that
%$X$ contains an isometric copy of $\{x,z \}*S_{z}$. In 
%particular, the antipodes of $x$ are discrete.
\end{proof}

\section{Local structure} \label{difsec}
In this section let $X$ be a fixed building-like space of dimension $n$.

\subsection{Simple remarks on the links} \label{simplerem}
 Let $x\in X$ be a point, $\gamma :[0,r] \to X$ a geodesic starting at $x$
in the direction $v\in S_x$.  From the local conicality of $X$ at $x$
we derive that for all small $t$ (for $t< \min \{ r, r_x \} $) the
link $S_{\gamma (t)} X$ has the
form $S_{\gamma (t)} X = S^0 * Z$, where $S^0 = \{ \gamma ^+ , \gamma ^- \}$
and $Z$ is canonically isometric to $S_v (S_x)$. 

 Assume now that $\gamma $ can be prolonged beyond $x$ to a geodesic
$\gamma :[-r,r] \to X$ and let $w\in S_x$ be the incoming direction of 
$\gamma$. Then $S_x$
contains an isometrically embedded $S_v(S_x) * \{ v,w \}$. Moreover
$S_v(S_x)= S_w (S_x)$. To see this observe that
$v$ and $w$ are antipodes in $S_x$. Each direction in $S_v (S_x)$ corresponds
to a direction in $S_{\gamma (t)} X$, orthogonal to $\gamma$.  Since 
$X$ is building-like at $\gamma (t)$, this direction gives a germ of spherical
triangle with one side $\gamma (-t)\gamma (t)$. This triangle defines
a geodesic in $S_x$ from $v$ to $w$ starting at the given direction of 
$S_v (S_x)$. Now the conclusion follows from the lune lemma 
(\cite[L.\ 2.5]{lunebal}).

\subsection{Regular points} We will call a point $x$ in $X$
 \emph{regular}, if $S_x =S^{n-1}$ holds. The set
of all regular points will be denoted by $R:=R(X)$. If~$x\in R$ is arbitrary,
 then $B_{r_x}(x)$ is isometric to a convex subset $C$ of the sphere $S^n$ and
 we have $S_x C =S^{n-1}$.  
In particular  $C$ contains a neighborhood
 of $x$ in $X$. Hence
the subset $R(X)$ is open in $X$ and locally
isometric to $S^n$.
 Due to Lemma \ref{closeness3} 
each convex subset of $R$ is a very spherical subset of $X$. 

We are going to prove that $R(X)$ is dense in $X$. The next lemma
is in fact true in arbitrary locally conical spaces.

\begin{lem}\label{sphereInDirections} 
Let $X$ be a building-like space with $\dim (X)=n$.
 Let $C\subset X$ be a convex subset with $\dim (C)= \dim (X)$. Then there
is a point $x\in C$ with $S_x C =S_x X = S^{n-1}$.
\end{lem}

\begin{proof}
 In dimension $n=1$ the statement follows directly from 
local conicality.
Let $n> 1$ be arbitrary. Choose $q\in C$ with $\dim (S_q C) =n-1$.
Due to \cite[Thm.\ B.3]{kleinerDimension} 
there must be a point $z\neq q $  with $z\in C \cap B_{r_q} (q)$ and
  $\dim (S_z C) =n-1$.
Then for each inner point $y$ of the geodesic $\gamma =qz$ we deduce
$\dim (S_y C) =n-1$, since there is a natural isometric embedding 
$S_z C\rightarrow S_{y}C$.
 However  $S_y C$ and $S_y X$ split as $S^0 * Z$ resp.
$S^0 * \tilde Z$, ($S^0 =\{ \gamma ^+ ,\gamma ^- \} $).
 Thus in a small
neighborhood of~$y$ the set $\tilde X$ (resp. $\tilde C$) 
of points   $\bar y \in X$ (resp. $\bar y \in C$)
with $ \angle_{y} (\bar y, q) =\frac \pi 2$ is a convex
subset of $X$ of dimension $n-1$. Arguing by  induction, 
 we find a
point $x \in \tilde C$ arbitrarily close to $y$  
with $S_x \tilde X = S_x \tilde C = S^{n-2}$.
The local structure of $X$ near $y$ implies  $S_x C = S_x X =S ^{n-1}$.
\end{proof}

We deduce

\begin{prop} \label{regexist}
Let $X$ be a building-like space  of dimension $n$. 
Then for each $x\in X$ holds $\dim (S_x ) =n-1$. Moreover $R(X)$ is dense
 in $X$.
\end{prop}

\begin{proof}
$R$ is open and non-empty by   \lref{sphereInDirections}.
Let $T$ be an open convex subset of regular points and $x\in X$. Due to
Lemma \ref{sphericalParts}  the convex hull $C$ of $x$ and $T$ is
spherical. We get 
$\dim (S_x) \geq \dim (S_x C ) \geq \dim (S) -1 = n-1$. 
Applying \lref{sphereInDirections} 
to an arbitrarily small neighborhood of $x$ we obtain
the second statement.
\end{proof}

\subsection{Easy applications} Let $X$ be a building-like space
of dimension~$n$ again. The remark at the end of Subsection \ref{basicbuild} 
 and the implication  $(4) \to (1)$ of \tref{introThm2} show
that, if $X$ is not a building,
then no regular point $x\in R(X)$ can have an antipode in $X$. 
The next result finishes the proof of \tref{introThm2}:

\begin{cor} \label{smallRadius}
 If $X$ is not a building, then
$ \rad (X) <\pi$.  
\end{cor}

\begin{proof}
Choose a regular point $x\in X$. For small $\varepsilon >0$ the ball 
$B_{\varepsilon} (x) $  consists of regular points. If for some
$z\in X$ we had $d(z,x) > \pi - \varepsilon$, then we could prolong the geodesic
$zx$ inside the ball $B_{\varepsilon} (x) $ (using \lref{easyLemma})
and obtain an antipode $\bar z$ of $z$ in $B_{\varepsilon} (x) $. Since 
$\bar z$ is regular, we would deduce that $X$ is a building.

 Thus if $X$ is not a building we obtain $\rad (X) \leq \pi  - \varepsilon$.
\end{proof}

The following result is mentioned in the introduction:
\begin{lem}\label{lowerDimSphere}
Let $X$ be as above.  Assume that $X$ contains a Euclidean sphere $S$
of dimension $n-1$.  Then either $X$ is a building or the radius of 
$X$ is $\frac \pi 2$.
\end{lem}

\begin{proof}
Let $x\in S$ be a point. Since $x$ has an antipode $y$ in $S$, we see
that $X$ contains   the isometrically embedded subset $S_x * \{x , y \}$.
Therefore  $S_x$ is an $(n-1)$-dimensional building-like space that contains
an $(n-2)$-dimensional sphere $S_x S$. Hence, an inductive argument shows 
  that $X$ contains an $n$-dimensional spherical hemi-sphere $H$
whose boundary sphere is $S$. Let  $m$ be the midpoint of this hemisphere.
Assume that $\rad _m (X) >\frac \pi 2$. Since the set of regular points is
 dense in $X$ we find a regular point $x$ with $d(x,m)> \frac \pi 2$. 
Consider the geodesic $xm$. By \lref{easyLemma} we
 find a direction $v \in S^{n-1} =S_m H$
that is antipodal to the starting direction of $mx$. Therefore we can 
prolong  the geodesic $xm$ inside $H$ and find an antipode of $x$ in $H$.
 This implies that $X$ is a building. 
\end{proof}

\subsection{Regular directions} Before we embark on the proof of 
\tref{boundarythm}
we will make some additional remarks about regular points. The picture
will be completed  in Subsection \ref{finalstr}.

\begin{defn}\label{regularDirection}
We call the starting  direction $v\in S_x$ of a geodesic $\gamma$
 \emph{regular} if for all sufficiently small $\varepsilon$ the point $\gamma
(\varepsilon )$ 
is regular.
\end{defn}

Remark that the direction $v=\gamma ^+\in S_x$ is regular 
iff $S_v (S_x) = S ^{n-2}$.
 Due to \pref{regexist}  the set of regular directions $R_x$ in $S_x$ is open
and dense. Moreover $R_x$ is locally isometric to $S^{n-1}$,
in particular it is locally compact.
 Using the second observation in Subsection \ref{simplerem}
 we see that if $x$ is an inner point
of the geodesic $\gamma$, then $\gamma ^+ \in S_x$ is regular iff the 
opposite direction $\gamma ^-$ is regular.

\subsection{Boundary}  In this subsection we are going to prove 
\tref{boundarythm}.
We start with a (technical) definition of the boundary:

\begin{defn} \label{auxdef}
A point $x$ in a building-like space $X$ will be called an
\emph{inner point} of $X$, if it is an inner point of a geodesic
connecting regular points. The set of all 
inner points will be denoted by $X_0$, its complement by $\partial X$.
\end{defn}

\begin{rem}  $X_0$ is in general not open in $X$, however it is open
with respect to the natural weak topology, see  \lref{weaktop}.
\end{rem}

\begin{lem}\label{sphereEquivalents}
Let $X$ be an $n$-dimensional building-like space.
Then  a point $x\in X$ is an inner point of $X$ iff there is
 an $n$-dimensional convex spherical subset $T$ that contains
$x$ as an inner point (of $T$).
\end{lem}

\begin{proof}
Assume that $x$ is an inner point of the geodesic $y_1y_2$ with regular
points $y_1y_2$. Then small neighborhoods $U_i$ of $y_i$ are 
very spherical. Hence $x$ is an inner point of the convex hull of $U_1$ 
and $U_2$ that is  spherical by \lref{sphericalParts}.

 On the other hand let $x$ be an inner point of $T$. By \pref{regexist} 
the set 
$R (X) \cap T$  of  regular points lying in $T$ is dense in $T$. Since it is 
also open and locally convex, one easily finds points $y_1,y_2 \in T\cap R$
such that $x$ lies on the geodesic $y_1y_2$.
\end{proof}

 Using \lref{easyLemma}  we immediately conclude:

\begin{cor}\label{oneEquivTwo}
A point $x$ is an inner point iff no geodesic terminates in $x$.
In particular  $X$ has
no boundary iff it is a building.
\end{cor}

Now let $x\in X$ be an inner point, let 
$\gamma :[0,\varepsilon' )\to B_{r_x} (x)$ be a geodesic starting at $x$.
Choose $T$ as in \lref{sphereEquivalents} and prolong $\gamma$ to
a geodesic $\gamma :[-r,\varepsilon' )$ inside $T$.  We see that a small
neighborhood~$V$ of $y=\gamma (-r)$ in $T$ and $\gamma (\varepsilon'' )$
(for $\varepsilon'' <\varepsilon '$) span
a spherical subset of~$X$. This shows that one can find an $n$-dimensional
 spherical convex subset~$T'$ containing $\gamma [0,\varepsilon ]$
(for $\varepsilon <\varepsilon '$).
Repeating this argument we see that for each geodesic
$\eta :(-\varepsilon' ,\varepsilon' )\to B_{r_x} (x)$ with $\eta (0)= x$
there is an $n$-dimensional convex spherical subset $T''$ that 
contains $\eta [-\varepsilon ,\varepsilon ]$.  
Now we are going to prove:

\begin{lem} \label{diffic}
A point $x$ is an inner point of $X$ iff $S_x$ is a building. 
\end{lem}

\begin{proof} Let $x$ be an inner point.
Denote by $V\subset S_x$ the set of all directions in which a geodesic
starts. 
By definition, $V$ is dense in $S_{x}$; furthermore, $V$ is convex
because of the locally conical structure of $X$. By the last observation,
 each pair of antipodes
$v,w\in V$ is contained in an $(n-1)$-dimensional sphere 
$S^{n-1} \subset V\subset S_x$.
 This implies that $V$ is geodesically complete (\lref{easyLemma}).

Unfortunately we cannot apply \tref{append}, since $V$  may not be  
complete.
To circumvent this difficulty, consider the ultraproduct 
$V^{\omega} =S_x ^{\omega}$. This is a (complete) CAT(1) space of 
dimension 
$n-1$  and each pair of antipodes of $V^{\omega}$ are still contained in
some $S^{n-1}$ (due to the geodesic completeness of $V$). Now $S_x$ contains
the (non-empty) open locally compact subset $R_x$ of regular points. 
Let $T\subset R_{x}$ be open and relatively compact. Then $T$ is also
an open and relatively compact subset of $V^{\omega}$.
%Hence
%$R_x$ is also an open relatively compact subset of $V^{\omega}$ and 
%applying
From \tref{append} (which is proved independently in section
\ref{veryend}), we deduce that $V^ {\omega}$ is a building. Hence $S_x$ 
is a convex subset of a building with 
$\rad (S_x) =\rad (S_x ^{\omega} ) =\pi$. Thus $S_x$ itself is a building
by \tref{introThm2}.  

 Assume now that $S_x$ is a building. Since each convex dense subset
in a building is the whole building, we see that in each direction 
$v\in S_x$ a geodesic starts. Choosing a finite number of directions 
in $S_x$ whose convex hull is a sphere of dimension $n-1$, we
deduce from local conicality that $x$ is contained in a spherical
$n$-dimensional subset. Hence $x$ is an inner point.
\end{proof}

\begin{proof}[Proof of \tref{boundarythm}]
$(1)\leftrightarrow (2)$ is Corollary \ref{oneEquivTwo}, and
$(1)\leftrightarrow (4)$ is the previous lemma.

Let $x\in X$ be a boundary point. 
Then some geodesic $yx$ cannot be extended beyond
$x$.
 Thus the contraction along geodesics starting at $y$ gives a 
contraction of $U\setminus \{ x \}$, where choosing $y$ 
arbitrarily close to~$x$ we may choose $U$ to be an arbitrarily small 
neighborhood of $x$.

 On the other hand if $x$ is an inner point of $X$, then the ball
$B_{r_x} (x)$ is a convex part of the building $S_x *S^0$. Remark 
that $B_{r_x} (x)$ contains  an $n$-dimensional
 spherical convex subset containing $x$. 
This subset defines a non-trivial element in the local
$(n-1)$-dimensional homology group at $x$ 
(compare \cite[sect.\ 6.2]{kleinerLeeb}). 
 This shows that 
$U\setminus \{ x \}$ cannot be contractible for small neighborhoods of $x$.

Finally $(5)\to (4)$ is clear and $(2)$ implies that the link $S_x$
has radius smaller than $\pi$ (as in the proof of 
\cref{smallRadius}), 
hence $(2)\to (5)$.
\end{proof}

%\begin{rem} \label{linktopo}
%In fact using \cref{oneEquivTwo} and arguing as in  
%\cref{smallRadius} one easiely deduces that the link $S_x$ has
%radius smaller than $\pi$ for  each boundary point $x\in X$ and is therefore
%contractible. Thus using \lref{diffic} we see
%that if $\dim (X) \geq 2$ then for each point $z\in X$ the link $S_z$ is
%connected.   
%\end{rem}

The next results   show that the boundary discussed above has the 
description announced in the introduction.
First, we show that $X_{0}$ is convex, in fact  even more is true:

\begin{lem} \label{lastconv}
 Let $x$ be an inner point of $X$, $\gamma:[0,s) \to X$ 
a geodesic starting at $x$. Then $\gamma$ consists of inner points of $X$. 
In particular $X_0$ is convex.
\end{lem}

\begin{proof}
The observation  preceding \lref{diffic} shows that
 the set $I$ of 
numbers $t$ with  $\gamma (t) \in X_0$ is open in $[0,s) $.
 Let $t \in (0,s)$ be a boundary point of $I$. Then for $\varepsilon
<r_{\gamma (t)}$, the point $\gamma (t-\varepsilon )$
is contained in a convex $n$-dimensional spherical subset $T$ and this
$T$ together with $\gamma (t+ \varepsilon )$ span an $n$-dimensional
spherical subset containing $\gamma (t)$. Thus~$t\in I$, in contradiction
to the openness of $I$.
\end{proof}

Since each convex dense subset of a building-like space
$X$ must contain all regular and therefore all inner points, we conclude

\begin{cor}
$\partial X$ is the
 largest subset of $X$  whose complement is convex and everywhere
 dense.
\end{cor}

%\begin{proof}
% Observe first that a convex dense subset $C$ of a building-like 
%space must contain all the regular points and therefore
%all inner points. 
%Since the regular points are dense, what remains to show is convexity
%of the interior. This follows directly from the next Lemma.
%\end{proof}

\subsection{Regular points revisited} \label{finalstr}
 Now we are going to investigate the combinatorial structure of $X$.
The next lemma shows that maximal very spherical subsets of $X$
define a decomposition of $X$ in cells.

\begin{lem}
Let $C$ be a maximal very spherical subset  of $X$. Denote by $C_0$
the set of inner points of $C$ as a spherical set. Then $C_0$ is open
in $X$. Moreover it is a connected component of $R (X)$.
\end{lem}

\begin{proof}
Let $m\in C_0$ be arbitrary. If $C_0$ does not contain a neighborhood of
$m$, then we can find a point $x $ close to $m$ that is not in $C_0$.
Due to Lemma \ref{closeness1}  the convex hull of $x$ and $C$ is a
very spherical subset of $X$  
in contradiction to the maximality of $C$. Since $C_0$ is open it is 
certainly contained in $R$. Therefore we only have to show that 
$C \setminus C_0 $ does not contain regular points. 
 Assume that $x$ is such a point. Then a neighborhood of $x$ is very spherical
and therefore we will find  points
 $m\in C_0$ and $x'\in X\setminus C$ that are close.  Using 
Lemma \ref{closeness1}  again we get a contradiction to the maximality
of $C$. 
\end{proof}

\begin{rem}
From \tref{boundarythm} it is easy to deduce that a point $x\in X$ is
a regular point iff it has a neighborhood homeomorphic to a manifold.
This shows that the decomposition in maximal very spherical subsets
is natural also from the topological point of view.
\end{rem}

\begin{lem} \label{weaktop}
 Let $C$ be a maximal very spherical subset of $X$. Then the intersection 
$X_0 \cap C$ is open in $C$. For $x\in X_0 \cap C$ some neighborhood
of $x$ in $C$  is isometric to an open subset in some Coxeter chamber. 
\end{lem}

\begin{proof}
Let $x$ be a point in $X_0 \cap C$. We know that $S_x$ is a building
and that in each direction $v\in S_x$ a geodesic starts. This shows that
$S_x C$ is a maximal chamber of the building $S_x X$. This chamber is
a convex hull of finitely many points, hence for a small number $\varepsilon $
in each direction $v\in S_xC$ a geodesic of length at least $\varepsilon$ starts.
This and \lref{lastconv} show that $X_0 \cap C$ is open in $C$ and
that the $\varepsilon$-ball around $x$ in $C$ is isometric to the
$\varepsilon $-ball around $x$ in the Coxeter 
chamber $S_x C * \{ x\} $.  
\end{proof}

\begin{rem}
 From this and well known facts about Coxeter groups in spheres it is easy to 
see that if a maximal very spherical subset $C$ of $X$ is contained in $X_0$,
then $C$ is a spherical join $C=C_1*C_2* \dotsb *C_k$, where each $C_i$ is 
either $1$-dimensional or isometric to the  Coxeter chamber of an 
irreducible Coxeter group.  
\end{rem}

 The following lemma is a weak equivalent of the statement that
the simplicial structure defined by the decomposition
in maximal very spherical subsets is thick. We leave the easy proof to the 
reader. 

\begin{lem}
Let $C$ be a maximal very spherical subset of $X$, $x\in C$ a point. If
$S_x C =\{ v  \} * S^{n-2}$ is a hemisphere, 
then either $x$ is a boundary point of
$X$ and  $C$ is a neighborhood of $x$, or
 $S_x X = T * S^{n-2}$, where the discrete  set $T$ has more than two points.
Moreover $x$ is contained in at least $3$ cells in the latter case. \qed
\end{lem}

\subsection{Maximal spherical subsets} This subsection is devoted to the 
proof of \pref{extend}.

We start with a criterion when a spherical subset is not maximal:

\begin{lem}
Let $C$ be a spherical subset of $X$, $m$ a point of
$C$, and $y$ be a point close to $m$. Let $v$ be the initial direction
of $my$. If the convex hull $H$ of $v$ and $S_{m}C$ is spherical, then
the convex hull of $C$ and $y$ is spherical.
\end{lem}

\begin{proof}
Observe that $C$ is mapped isometrically by $p_{m}:X\rightarrow
S_{m}*S^{0}$.

Since every triangle $myq$ for $q\in X$ is spherical, we
find that $\{y \}\cup C$ is mapped isometrically (by $p_{m}$) into the
spherical set 
$H*S^{0}$. This implies that for all
$x_{1},x_{2}\in C$, the triangle $yx_{1}x_{2}$ is spherical.
Hence, the convex hull of $C$ and $y$ is spherical by
\lref{sphericalParts}. 
\end{proof}

\begin{proof}[Proof of \pref{extend}]
Since the closure and a union of a chain of spherical subsets is spherical,
we may assume that the spherical subset~$C$ of $X$ is maximal; we have 
to prove that it has dimension $n$ and that $\partial C \subset \partial X$.

 Let $m \in \partial C$ be an inner point of $X$.  Then $S_m X$ is a
building, hence (by induction or 
due to \cite[Prop.\ 3.9.1]{kleinerLeeb}) $S_m C\subsetneq
S^{n-1}\subset S_{p}X$, which is a contradiction to the maximality of
$C$ by the previous lemma  (since for an inner point of $X$ a
geodesic starts in all directions).

 Assume now that $C$ has dimension smaller than $n$ 
(the same argument as above implies that $C$ must be contained in
$\partial X$ in  
this case).   Pick an inner point $m$ of $C$. 
Similar to the proof of Lemma \ref{sphereInDirections}, one shows that
one can find another inner point $\bar m$ of  $C$, such that
$S_{\bar m} C$ is a (non-trivial) spherical join factor of $S_{\bar m} X$. 
This shows the existence of  a vector $v$ as in the lemma 
above and therefore a contradiction to the maximality.
\end{proof}

\section{Appendix} \label{veryend}
Here, we are going to prove \tref{append} and \cref{kleiner}:

\begin{proof}[Proof of \tref{append}]
By the assumptions $X$ contains at least one isometrically embedded $S^n$.
Due to Lemma \ref{easyLemma}  each point has an antipode and is
therefore contained in a  
sphere of maximal dimension. Hence, $X$ is geodesically
complete and has diameter $\pi$. 

We proceed by induction on dimension. In dimension $0$ there is nothing
to be done. If $X$ is reducible, i.e. if it has a non-trivial decomposition
$X=Y*Z$, then one easily sees, that the assumptions of the theorem are
fulfilled for the spaces $Y$ and $Z$, which have dimension smaller than $X$.
By induction we obtain that $Y$, $Z$ and therefore
$X=Y*Z$ are buildings. Hence we may assume that $X$ is irreducible.

Let $U$ be an open relatively compact subset of $X$.
Since $X$ is geodesically complete, $U$ contains an open subset $\tilde U$
homeomorphic to a manifold (compare for example \cite{Otsu}).
 The dimension of $\tilde U$ is at most 
$\dim (X) =n$ and since each point of $\tilde U$ is contained in some $S^n$,
we see that $\dim (\tilde U) =n$ and that $\tilde U$ is locally isometric
to  $S^n$.

 Consider the set $O$ of all points $x\in X$ that have a neighborhood
isometric to an open subset of $S^n$. This set is open by definition and
we have just seen that it is non-empty.
 If $O$ is the whole set $X$ then $X$ is 
isometric to the sphere $S^n$, hence we may assume $O\neq X$.

  Let $x\in O$ be arbitrary, $y\in X$ a point with $d(x,y)=\pi$.
Let $S=S^n$ be a sphere of maximal dimension that contains $x$ and $y$.
Since $x\in O$, the sphere $S$ contains a ball  $B_r (x) \subset
O\subset X$ for some 
small $r>0$. 
We are going to show that $S$ contains  $B_r (y)$.

Let $\bar y\in B_{r}(y)$.
Due to Lemma \ref{easyLemma}  we may continue
the geodesic $\bar yy$ inside $S$ and obtain an antipode
$\bar x$ of $\bar y$ in $S$. We have $d(x,\bar x) =d(\bar y, y)\leq r$,
hence  $\bar x \in O$.  
So there is a spherical neighborhood $U$ of $\bar{x}$. 
Since $S$ is  a sphere in $X$ containing $\bar{x}$, we have $U\subset
S$. 

If we assume $U$ to be a maximal connected spherical neighborhood of
$\bar{x}$, we have $U\supset 
B_{r}(x)$. Hence, the geodesic segment $\bar{x}x$ can be extended to a
geodesic $\bar{x}x\bar{y}$, implying that $x\bar{y}y$ is a geodesic
too. In particular, $\bar{y}\in S$.

Thus $S$ contains a neighborhood of $y$, and $y$ is in $O$. This shows 
that  the complement $T=X\setminus O$ is closed
and contains all antipodes of all of its points. Since $T\neq X$ and $X$
is irreducible we may apply  the main result of \cite{lytchakSphB}, which 
says that the existence of such a subset in an irreducible geodesically
complete space $X$  implies that $X$ is   a building.
\end{proof}

\bigbreak

In the proof of \cref{kleiner}, we use the appropriate definition of
local conicality in the $CAT(0)$ setting.

\begin{proof}[Proof of \cref{kleiner}]
 \lref{easyLemma} shows that $X$ is geodesically
complete. Since $X$ is proper, the sequence $(tX,x)$ converges
for $t\to \infty$
to the tangent cone $CS_x$ in the Gromov-Hausdorff topology and from 
\tref{append} we immediately obtain that each link is a spherical building.
By \cite{charneyLytchak} it is enough to prove that $X$  is locally
conical. 

Let $x\in X$ be arbitrary. Then for some $r>0$ each $n$-dimensional
Euclidean $\R ^n =F\subset X$ with $d(x,F)<r$ must contain $x$, since 
otherwise we would obtain a flat $\bar F =\R^n$ in the Euclidean cone $CS_x$
that does not contain the origin and this would contradict  
$\dim (CS_x) \leq n$. 

Hence for all $y,z\in B_r (x)$ each maximal flat through $y$ and $z$ must
contain $x$, thus the triangle $xyz$ is flat. Therefore   $X$ is locally 
conical.
\end{proof}

\bibliographystyle{alpha}
\bibliography{build}

\end{document}